\documentclass[12pt, Russian]{article}
\usepackage[russian]{babel}
\evensidemargin=\oddsidemargin \textwidth=165mm \textheight=238mm
\usepackage[dvips]{graphicx}
\usepackage{latexsym}
\usepackage{amssymb}
\begin{document}
\null
\begin{flushleft}\bf{\small}
\end{flushleft}
\vskip 5mm

\begin{center}
    \Large \bf Estimates of solutions to the linear Navier-Stokes equation
       \end{center}
\begin{center}
        {\sl Dr. Bazarbekov Argyngazy B.     \\}
        {\it Kazakh National State University, EKSU.~ e-mail : arg.baz50@mail.ru}

\end{center}

\textbf{Abstact~.}\begin{scriptsize} The linear Navier-Stokes
equations in three dimensions are given by: $u_{it}(x,t) - \rho
\triangle u_i(x,t) - p_{x_i}(x,t) = $ ~~~ $=w_i(x,t)$ ,
$div~\textbf{u}(x,t) = 0~~ ,~i = 1,2,3$~ with initial conditions:
$\textbf{u}|_{(t=0)\bigcup
\partial\Omega} = 0$. The Green function to the Dirichlet problem $\textbf{u}|_{(t=0)\bigcup
\partial\Omega} = 0$ of the equation  $u_{it}(x,t) - \rho
\triangle u_i(x,t)  = f_i(x,t)$ present as: $G(x,t;\xi,\tau) =
Z(x,t;\xi,\tau)~+~V(x,t;\xi,\tau).$~Where $Z(x,t;\xi,\tau) =
\frac{1}{8 \pi^{3/2} (t - \tau)^{3/2}}\cdot e^{- \frac{(x_1 -
\xi_1)^2 + (x_2 - \xi_2)^2 + (x_3 - \xi_3)^2 }{4 (t - \tau)}}$ - is
the fundamental solution to this equation [1 p.137] and
$V(x,t;\xi,\tau)$ -is the smooth function of variables
$(x,t;\xi,\tau)$. The construction of the function $G(x,t;\xi,\tau)$
is resulted in the book  [1 p.106]. By the Green function we present
the Navier-Stokes equation as: ~~$u_i(x,t)=\int_0^t
\int_{\Omega}\Big(Z(x,t;\xi,\tau) + V(x,t;\xi,\tau)\Big) \frac{d
p(\xi,\tau)}{d \xi} d \xi d \tau + \int_0^t
\int_{\Omega}G(x,t;\xi,\tau)w_i(\xi,\tau) d \xi d \tau$. But  $div
\textbf{u}(x,t)=\sum_1^3 \frac{d u_i(x,t)}{dx_i}=0.$ ~Using these
equations and the following properties of the fundamental function:
$Z(x,t;\xi,\tau)$: ~$\frac{d Z(x,t;\xi,\tau)}{d x_i} = - \frac{d
Z(x,t; \xi,\tau)}{d \xi_i},$~for the definition of the unknown
pressure  p(x,t) we shall receive the integral equation. From this
integral equation we define the explicit  expression of the
pressure: $p(x,t)~= -\frac{d }{d t} \triangle^{-1}\ast \int_0^t
\int_{\Omega}\sum_1^3 \frac{dG(x,t;\xi,\tau)}{d x_i}w_i(\xi,\tau) d
\xi d \tau + \rho \cdot \int_0^t \int_{\Omega}\sum_1^3
\frac{dG(x,t;\xi,\tau)}{d x_i}w_i(\xi,\tau) d \xi d \tau.$ By this
formula the following
estimate:~$\int_0^t\sum_1^3\Big\|\frac{\partial p(x,\tau)}{\partial
x_i}\Big\|_{L_2(\Omega)}^2 d \tau~
 <~ c \cdot \int_0^t\sum_1^3\|w_i(x,\tau)\|_{L_2(\Omega)}^2 d \tau
$ holds.

\end{scriptsize}

                                               \label {Baizhuman}

\vskip 4mm

\begin{scriptsize} 2000 Mathematical Subject Classification. Primary:~$35~K~55$ ;
Secondary:  $46 E 35$;  Keywords: postulate , integral equation,
Green function, Dirichlet problem, operator , smooth functions.
\end{scriptsize}

\textbf{1.~Introduction}.~Let $\Omega \subset R^3$ be a finite
domain bounded by the Lipschitz surface $\partial \Omega$.~$Q_t =
\Omega \times [0,t], ~x = (x_1,x_2, x_3)$ and
$\textbf{u}(x,t)~=~(u_i(x,t)_{i=1,2,3}$ are the vector functions.
~~The boundary problems for linear Navier-Stokes equations were
considered in the book [2 p.89-96].~If $\partial \Omega$ -is a
smooth surface, then there exists a solution: $\textbf{u}(x,t):
\textbf{u}_{x_i x_j}(x,t) , p(x,t) \in \textbf{L}_2(\Omega')\times
[0,T]$~where~$\Omega'\subset \Omega$~is a subarea of the area
$\Omega$. ~~I do not know other works devoted to this problem.

\vskip 1mm \textbf{2.~~The formulation and proof of result. }

Let~$\Omega \subset R^3$~be a finite domain bounded by the Lipschitz
surface ~$\partial\Omega.$~$Q_{t} = \Omega \times[0,t] , ~x
=(x_1,x_2,x_3)$~and ~$\textbf{u}(x,t)~=~ (u_i(x,t)_{i = 1,2,3}~,$
$\textbf{f}(x,t) = (f_i(x,t)_{i =1,2,3}$- are vector functions.
Here~t > 0, $\rho > 0$ are an arbitrary real numbers.~ The linear
Navier-Stokes equations are given by:

$$\frac{\partial u_i(x,t)}{\partial t}
~-~\rho~\triangle~u_i(x,t)~-~\frac{\partial p(x,t)}{\partial
x_i}~~=~w_i(x,t)~    \eqno(2,1)$$
$$div~\textbf{u}(x,t)~=~\sum_{i=1}^3\frac{\partial u_i(x,t)}{\partial x_i}~=~
0~~,i~=~1,2,3     \eqno(2,1')$$

$$\textbf{u} (x,0) =~0~~, ~~~ \textbf{u}(x,t)\mid_{\partial \Omega \times [0,t]}~=~0 \eqno(2,1")$$
The Green function for Dirichlet problem is as follows:
$$G(x,t;\xi,\tau) = Z(x,t;\xi,\tau)~+~V(x,t;\xi,\tau)$$~where
$Z(x,t;\xi,\tau) = \frac{1}{8 \pi^{3/2} (t - \tau)^{3/2}}\cdot e^{-
\frac{(x_1 - \xi_1)^2 + (x_2 - \xi_2)^2 + (x_3 - \xi_3)^2 }{4 (t -
\tau)}}$ ~is the fundamental solution for this equation  [1 p.106]
and  $V(x,t;\xi,\tau)$ is the  \textbf{smooth} function of variables
$(x,t;\xi,\tau)$. The construction of the function $V(x,t;\xi,\tau)$
is resulted in the book  [1 p.106]. By the Green function we present
the Navier-Stokes equation as: ~~
$$u_i(x,t)=\int_0^t\int_{\Omega}\Big(Z(x,t;\xi,\tau) + V(x,t;\xi,\tau)\Big) \frac{d
p(\xi,\tau)}{d \xi} d \xi d \tau + \int_0^t
\int_{\Omega}G(x,t;\xi,\tau)w_i(\xi,\tau) d \xi d \tau$$.
$$i = 1,2,3.~~~~~~~~~~~~~~~div \textbf{u}(x,t)=\sum_1^3 \frac{d u_i(x,t)}{dx_i}=0
\eqno(2,2)$$ ~~We shall prove the following fundamental result:

\textbf{Theorem  2,1} There exists the explicit expression to the
pressure p(x,t),  depending  on  the right-hand side $w_i(x,t)$:
$$p(x,t)~=-\frac{d }{d t}
\triangle^{-1}\ast \int_0^t \int_{\Omega}\sum_1^3
\frac{dG(x,t;\xi,\tau)}{d x_i}w_i(\xi,\tau) d \xi d \tau + \rho
\cdot \int_0^t \int_{\Omega}\sum_1^3 \frac{dG(x,t;\xi,\tau)}{d
x_i}w_i(\xi,\tau) d \xi d \tau$$where  $\triangle^{-1}$ - is the
inverse operator to the Dirichlet problem of Laplase  operator for
the domain $\Omega$.~~ To the solutions  $u_i(x,t), ~p(x,t)$ for
these linear Navier-Stokes equation the following estimates:
$$\|\textbf{u}(x,t)\|_{W_2^{2,1}(Q_t)}~<~c \cdot
\|\textbf{w}(x,t)\|_{L_2(Q_t)}$$
$$\int_0^t\sum_1^3\Big\|\frac{\partial p(x,\tau)}{\partial x_i}\Big\|_{L_2(\Omega)}^2 d \tau~
 <~ c \cdot \int_0^t\sum_1^3\|w_i(x,\tau)\|_{L_2(\Omega)}^2 d \tau
\eqno(2,3)$$ are valid.~~Here and below by symbol  $c$~we denote a
generic constants, independent on the solution and right-hand side
whose value is inessential to our  aims, and it may change from line
to line. $\blacktriangleleft$

\textbf{Postulate:} Let   t > 0 be an arbitrary real number. ~For
the proof of the estimates  (2,3) we assume that: ~$p(x,t) ~\in
~L_2(W_2^1(\Omega);[0,t]).$~~~I.e.
$$\int_0^t\sum_1^3\Big\|\frac{\partial p(x,\tau)}{\partial x_i}\Big\|_{L_2(\Omega)}^2 d \tau~
 <~\infty \eqno(2,4)$$\textbf{Remark 1.} For the function $p(x,t) \in
~L_2(W_2^1(\Omega);[0,t])$ there exists~ $p(x,t)|_{\partial \Omega}
\in L_2(\partial \Omega;[0,t]).$$\blacktriangleleft$

 Below, for the definition of the unknown pressure
 $p(x,t): p(x,t) \in L_2(Q_t)$ we shall receive the integral equation. And from this integral equation
 we shall receive the explicit  expression of the pressure p(x,t) and estimates (2,3).

It is known , that $G(x,t;\xi,\tau)=G(\xi,t;x,\tau)$ [3 p. 305] and
 $G(x,t;\xi,\tau)|_{\xi \in \partial\Omega} = G(x,t;\xi,\tau)|_{x
\in \partial\Omega} = 0.$~ Since $G(x,t;\xi,\tau)|_{\xi \in
\partial\Omega} = 0$ and   $p(x,t) \in
~L_2(W_2^1(\Omega);[0,t]),$~ integrating by part, from  (2,2) we
have:
 $$u_i(x,t)~=~- \int_0^t \int_{\Omega}\Big(\frac{d Z(x,t;\xi,\tau)}{d \xi_i} +
\frac{d V(x,t;\xi,\tau)}{d \xi_i}\Big)  p(\xi,\tau)d \xi d \tau
+$$$$+ \int_0^t \int_{\Omega}G(x,t;\xi,\tau)w_i(\xi,\tau) d \xi d
\tau \eqno(2,5)$$ Since  $\frac{d Z(x,t;\xi,\tau)}{d x_i} = -
\frac{d Z(x,t; \xi,\tau)}{d \xi_i}$ and $div
\textbf{u}(x,t)=\sum_1^3 \frac{d u_i(x,t)}{dx_i}=0$, from   (2,5) we
obtain:
$$\triangle_x\int_0^t\int_{\Omega}Z(x,t;\xi,\tau)p(\xi,\tau) d\xi
d\tau~-~\int_0^t\int_{\Omega}\Big(\sum_1^3\frac{d^2V(x,t;\xi,\tau}{dx_id\xi_i}\Big)
p(\xi,\tau) d \xi d \tau +$$ $$+\int_0^t \int_{\Omega}\sum_1^3
\frac{dG(x,t;\xi,\tau)}{d x_i}w_i(\xi,\tau) d \xi d \tau~=~0
\eqno(2,6)$$Since  $\sum_1^3\frac{d}{d x_i} \frac{d
V(x,t;\xi,\tau)}{d x_i}~=~ \bigtriangleup_x V(x,t;\xi,\tau)$, we
rewrite this equation in the following form:

$$\triangle_x\int_0^t\int_{\Omega}\Big(Z(x,t;\xi,\tau) + V(x,t;\xi,\tau)\Big) p(\xi,\tau) d\xi
d\tau~- $$$$-\int_0^t \int_{\Omega} \sum_1^3 \frac{d}{d
x_i}\Big(\frac{d V(x,t;\xi,\tau)}{d \xi_i} +
\frac{dV(x,t;\xi,\tau)}{d x_i}\Big) p(\xi,\tau) d \xi d \tau =
\eqno(2,7)
$$$$ - \int_0^t \int_{\Omega}\sum_1^3 \frac{dG(x,t;\xi,\tau)}{d
x_i}w_i(\xi,\tau) d \xi d \tau $$ Let us denote:
$$\bigtriangleup \ast u(x,t) =~\sum_1^3 u_{x_i x_i}(x,t) ;~~~~~~~;\bigtriangleup^{-1} \ast f(x,t)~=~\int_{\Omega}
G_D(x,\xi)~f(\xi,t) d \xi$$
$$~T\ast u(x,t) = \frac{\partial u(x,t)}{\partial t}
-\rho\cdot \triangle~u(x,t);~~~~~~~~~; T^{\ast}\ast u(x,t) =
 -\frac{\partial u(x,t)}{\partial t} - \rho~\triangle~u(x,t)
\eqno(2,8)$$
$$V \ast p(x,t)~=~\int_0^t\int_{\Omega}\sum_1^3\frac{d}{dx_i}\Big(\frac{dV(x,t;\xi,\tau)}{d\xi_i}
+ \frac{dV(x,t;\xi,\tau)}{d x_i}\Big) p(\xi,\tau) d \xi d \tau $$
where ~$G_D(x;\xi)$~is the Green function of the Dirichlet problem
for Laplase operator and the domain $\Omega$.~~~ Note that:
~$\int_0^t\int_{\Omega}\Big(Z(x,t;\xi,\tau) + V(x,t;\xi,\tau)\Big)
p(\xi,\tau) d\xi d\tau \Big|_{x \in \partial \Omega}~=~0$.~~Using
these denotes , we rewrite the basis equation (2,7):
$$\int_0^t\int_{\Omega}\Big(Z(x,t;\xi,\tau) + V(x,t;\xi,\tau)\Big)
p(\xi,\tau) d\xi d\tau~-~\bigtriangleup^{-1}\ast~V\ast
 p(x,t)~=$$$$=~-\bigtriangleup^{-1}\ast\int_0^t \int_{\Omega}\sum_1^3
\frac{dG(x,t;\xi,\tau)}{d x_i}w_i(\xi,\tau) d \xi d \tau
\eqno(2,9)$$Since ~$T\ast \int_0^t\int_{\Omega}\Big(Z(x,t;\xi,\tau)
+ V(x,t;\xi,\tau)\Big) p(\xi,\tau) d\xi d\tau~=~p(x,t)$, from (2,9)
we have:
$$p(x,t) - T\ast  \bigtriangleup^{-1}\ast V\ast
 p(x,t) = - T\ast\bigtriangleup^{-1}\ast\int_0^t \int_{\Omega}\sum_1^3
\frac{dG(x,t;\xi,\tau)}{d x_i}w_i(\xi,\tau) d \xi d \tau
\eqno(2,10)$$ \textbf{Lemma  2,1.}~~For any function  $p(x,t) \in
L_2(Q_t)$ ~ the following equality is valid:
$$T_{(x,t)} \ast V\ast
 p(x,t)~=~\frac{d~ V\ast p(x,t)}{d t} ~-~\rho\cdot \triangle_x V\ast
 p(x,t)  ~\equiv~0   \eqno(2,11)$$
where the operator $V\ast p(x,t)$ is defined above by the formula:
 $$V \ast p(x,t)~=~\int_0^t\int_{\Omega}\Big(\sum_1^3
\frac{d^2V(x,t;\xi,\tau)}{dx_i d\xi_i} + \bigtriangleup_x
V(x,t;\xi,\tau)\Big) p(\xi,\tau) d \xi d \tau   \eqno(2,12)$$
$\blacktriangleright$ On the domain  $\Omega \times [0,t]$
А.Friedman has constructed in the book  [1 p. 106-111] the Green
function of Dirichlet problem to the parabolic equation: ~$$T \ast
u(x,t) = u_t(x,t) - \rho \cdot \bigtriangleup_x u(x,t) =
p(x,t);~~u_{\{t = 0\}\bigcup
\partial\Omega \times [0,t]}~=~0.$$~And present the solution to this problem as:$$u(x,t)=\int_0^t
\int_{\Omega} \Big(Z(x,t;\xi,\tau) + V(x,t;\xi,\tau)\Big)p(\xi,\tau)
d \xi d\tau    \eqno(2,13)$$where the function  $V(x,t;\xi,\tau)$ by
variables   $(\xi,\tau)$  satisfies to the following equation and
boundary conditions [1 p. 108]:
$$T^{\ast}_{(\xi,\tau)} \ast V(x,t;\xi,\tau)~=~- \frac{d}{d
\tau}V(x,t;\xi,\tau) ~-~\rho \cdot \triangle _{\xi}V(x,t;\xi,\tau)
~=~0  \eqno(2,14)$$
$$V(x,t;\xi,\tau)\Big|_{(\xi,\tau) \in \partial \Omega \times [0,t]}
~=~-~Z(x,t;\xi,\tau)\Big|_{(\xi,\tau) \in \partial \Omega \times
[0,t]}   \eqno(2,14')$$$$V(x,t;\xi,\tau)\Big|_{\tau = t}~=~V(x,t;
\xi,t)~=~0   \eqno(2,14")$$ The formula (2,13) follows from
integrating by part the following expression:
$$\int_0^t\int_{\Omega}\Big( Z(x,t;\xi,\tau) +
V(x,t;\xi,\tau)\Big)\Big(u_{\tau} ~-~\rho\cdot
\triangle_{\xi}u(\xi,\tau)\Big) d \xi d \tau ~=~u(x,t),$$ using the
boundary conditions: (2,1"), (2,14') , (2,14").

The function  $V(x,t;\xi,\tau)$  by variables  $(x,t)$ satisfies to
the following parabolic equation: $T_{(x,t)}\ast
 V(x,t;\xi,\tau)~ =~0$.~Then, the following equalities are valid: $$T_{(x,t)}\ast
 \frac{d^2 V(x,t;\xi,\tau)}{d x_i d \xi_i} = \frac{d}{d x_i}
 T_{(x,t)}\ast \frac{d V(x,t;\xi,\tau)}{d \xi_i}= 0$$$$
 T_{(x,t)}\ast\triangle_x V(x,t;\xi,\tau) = \triangle_x \ast T_{(x,t)}\ast V(x,t;\xi,\tau)~
=~0$$Since~$V(x,t;\xi,t) = 0$~ for any $x \in \Omega$,~$\xi \in
\Omega$, then~$\frac{d^2 V(x,t;\xi,t)}{dx_i d \xi_i} = 0$ and
~$\triangle_x V(x,t;\xi,t) = 0.$  Therefore, from (2,12) we have
~$T_{(x,t)} \ast V\ast
 p(x,t) = 0.$ Lemma 2,1 is proved.$\blacktriangleleft$

\textbf{The Laplase and parabolic opеrаtоrs.} The Dirichlet problem
for Laplase operator has infinite set of the own functions $u_n(x)$
and own numbers $\lambda_n > 0:$   $\triangle u_n(x)=-\lambda_n\cdot
u_n(x); ~~ |u_n(x)|_{L_2(\Omega)}=1$~,~
 $(u_n(x),u_m(x))_{L_2(\Omega)} = 0.$~ The each solution  u(x,t)
to the homogeneous parabolic equation  $T\ast u(x,t) = 0$ present
as:$$u(x,t)~=~\sum_1^{\infty} c_n\cdot e^{- \lambda_n \cdot
\rho\cdot t}\cdot u_n(x) \eqno(2,15)$$It follows from the
decomposition of the function  u(x,t) on the own functions
$\{u_n(x)\}.$

\textbf{Lemma 2,2.}~For any function  $p(x,t) \in L_2(Q_t)$ ~is
valid the following equality:
$$T_{(x,t)} \ast \bigtriangleup^{-1} \ast V\ast
 p(x,t) ~~\equiv~~0   \eqno(2,16)$$
$\blacktriangleright$ Since ~$\triangle \ast \triangle^{-1}\ast
f(x,t) ~=~\triangle_x \ast\int_{\Omega}G_D(x,\xi)f(\xi,t) d \xi
~=~f(x,t),$~ it follows by the definition of the operators
$T,\bigtriangleup^{-1}$ that:
$$T\ast\bigtriangleup^{-1} \ast V\ast p(x,t)~=~\frac{d }{d t}
\triangle^{-1}\ast V\ast p(x,t) ~-~\rho \cdot V\ast p(x,t)
\eqno(2,17)$$ As it is proved in Lemma 2.1:~ $T\ast V\ast
p(x,t)~=~0.$ Decomposing the function  $V\ast p(x,t)$ on the own
functions to Laplase operator, we shall receive his representation
as:
$$V \ast p(x,t)~=~\sum_1^{\infty} c_n\cdot e^{- \lambda_n \cdot
\rho\cdot t}\cdot u_n(x)   \eqno(2,15')$$ Since $\triangle^{-1}
u_n(x) = - \frac{1}{\lambda_n}\cdot u_n(x)$, then
$\bigtriangleup^{-1}\ast V\ast p(x,t) = - \sum_1^{\infty}
\frac{c_n}{\lambda_n}\cdot e^{-\rho\cdot \lambda_n \cdot t} \cdot
u_n(x).$ ~Substituting this function to the formula (2,17), we shall
receive that:
 $T_{(x,t)} \ast \bigtriangleup^{-1} \ast V\ast
 p(x,t) ~~\equiv~~0.$~~Lemma  2.2 is proved.$\blacktriangleleft$

From the equation (2,10) , Lemma 2.2 and formula  (2,17) we shall
receive the following explicit expression to pressure function
p(x,t), depending  on the right-hand side \textbf{w}(x,t):
$$p(x,t) ~ =~ - T\ast\bigtriangleup^{-1}\ast\int_0^t \int_{\Omega}\sum_1^3
\frac{dG(x,t;\xi,\tau)}{d x_i}w_i(\xi,\tau) d \xi d \tau ~=
\eqno(2,18)$$
$$=-\frac{d }{d t}
\triangle^{-1}\ast \int_0^t \int_{\Omega}\sum_1^3
\frac{dG(x,t;\xi,\tau)}{d x_i}w_i(\xi,\tau) d \xi d \tau + \rho
\cdot \int_0^t \int_{\Omega}\sum_1^3 \frac{dG(x,t;\xi,\tau)}{d
x_i}w_i(\xi,\tau) d \xi d \tau$$Let us notice that by this formula
follows the following formula: $\triangle_x p(x,t) = \sum_1^3
\frac{w_i(x)}{d x_i}.$ The same formula follows from the
Navier-Stokes equation.~~ It is obvious , that:
$$\int_0^t \int_{\Omega}\sum_1^3 \frac{dG(x,t;\xi,\tau)}{d
x_i}w_i(\xi,\tau) d \xi d \tau~\in~~W^{1,1}_2(Q_t)$$ It follows by
the second formula (2,18) the following estimate:
$$\int_0^t\Big(\Big\|T\ast\bigtriangleup^{-1}\ast\int_0^{\tau}
\int_{\Omega}\sum_1^3 \frac{d G(x,t;\xi,\tau_1)}{d
x_i}w_i(\xi,\tau_1) d \xi d \tau_1\Big\|_{W_2^1(\Omega)}\Big)^2 d
\tau ~<~c~|\textbf{w}(x,\tau)|^2_{L_2(Q_t)}$$ Differentiating the
equation  (2,18) by $x_i$, from this estimate we obtain an estimate
(2,3).

\textbf{Remark 2.} By the formula  (2,18) we find  $\frac{d
p(x,t)}{d x_i}, i = 1,2,3$, depending on the functions $w_i(x,t)$,
and substitute these functions to the Navier-Stokes equation (2,2):
~ $u_i(x,t)=\int_0^t\int_{\Omega}\Big(Z(x,t;\xi,\tau) +
V(x,t;\xi,\tau)\Big) \frac{d p(\xi,\tau)}{d \xi} d \xi d \tau +
\int_0^t \int_{\Omega}G(x,t;\xi,\tau)w_i(\xi,\tau) d \xi d \tau$. We
find the explicit expression to the solutions $u_i(x,t)$ by the
right-hand side $w_i(x,t)$.~~ Theorem (2,1) is
proved.$\blacktriangleleft$~~ I am sorry for taking you time.

\newpage

\begin{center}  R~e~f~e~r~e~n~c~e~s.\end{center}

[1]. Friedman Avner. : Partial differential equations of parabolic
type.,

~~~~~Prentice-Hall, (1964).

[2]. Ladyzhenskaya O.: The mathematical theory of Viscous
 Incompressible Flows,

 ~~~~~Cordon and Breach , ~(1969).

[3].Sobolev S.L.Partial differential equations of Mathematical
Physics.,Moskow,(1966).

\vskip35mm

Address:050026 ,Muratbaeva St. 94-5,

Almaty , Republic Kazakhstan.\label{BAZARBFIN}

\end{document}